\documentclass[11pt]{article}
\usepackage{latexsym}
\usepackage[dvips]{epsfig}

\newcommand{\qed}{\mbox{$\Diamond$}\vspace{\baselineskip}}

\newtheorem{theorem}{Theorem}[section]

\newtheorem{definition}[theorem]{Definition}

\newenvironment{proof}{\noindent {\bf Proof:}}{{\qed}}

\newcommand{\vanish}[1]{}
\begin{document}

\title{An (almost) optimal answer to a question by Wilf}
\author{Micah Coleman \thanks{
        Department of Mathematics,
        University of Florida,
        Gainesville FL 32611-8105. 
        Supported by a grant from the 
	University of Florida
	Undergraduate Scholars Program, 
	mentored by Mikl\'os B\'ona.}}
\maketitle

\date{}
\begin{abstract}
We present a class of permutations for which the number 
of distinctly ordered subsequences of each permutation approaches an almost 
optimal value as the length of the permutation grows to infinity. 
\end{abstract}

\section{Introduction}

\begin{definition}
We say that two permutations $q = q_1\ q_2\ \cdots\ q_k$ and $r = r_1\ r_2\ \cdots\ r_k$ are {\em identically ordered} provided $q_i < q_j \iff r_i < r_j$
\end{definition}

\begin{definition}
Let $p$ be a permutation of length $n$ with entries denoted 
$p_1, p_2, \ldots, p_n$. 
We define a {\em subsequence} of $p$ to be a sequence 
$p_{i_1}\  p_{i_2}\ \ldots\ p_{i_k}$ where $i_1 < i_2 < \cdots < i_k$.
We say that $p\ $ {\em contains} a {\em pattern} $q = q_1\ q_2\ \cdots\ q_k$ if
there exists a $k$-subsequence of $p\ $ that is ordered identically with $q$.
\end{definition}

At the Conference on Permutation Patterns, Otago, 
New Zealand, 2003, Herb Wilf asked how many distinct patterns could there be in a permutation of length $n$.  Based on empirical evidence, it seems this number may approach the theoretical upper bound of $2^n$.  In this paper we enumerate patterns contained in a certain class of permutations to at least establish a lower bound for this function.

Let $f(p)$ be the number of distinct patterns contained in a  
permutation $p$. 
Let $h(n)$ be the maximum $f(p)$ for any permutation $p$ of length $n$.
Upper bounds for $h(n)$ are the number of all subsequences, 
$2^n$, and the number of permutations for each pattern
length, $k!$. 
As $n$ grows, this second bound quickly becomes insignificant, as
\[{n \choose k} < k!\] for all $k$ above a breakpoint which grows
much slower than $n$.  Therefore, we are concerned with proving
that $h(n)$ grows almost as fast as $2^n$.

Wilf demonstrated a class of permutations for which $f > ({1 + \sqrt{5} \over 2})^n$.

Let $W_n$ denote the $n$-permutation 
$(1\ \ n\ \ 2\ \ n-1\ \ \cdots\ \ \lfloor {n \over 2}+1 \rfloor)$.
It should be clear that 
$$f(p_1\ \ p_2\ \ \cdots\ \ p_n) = 
f(n-p_1+1\ \ n-p_2+1\ \ \cdots\ \ n-p_n+1)$$
Then, the number of patterns contained in $W_n$ which do not include the 
first entry is $f(W_{n-1})$. 
The number of patterns contained in $W_n$ which are required to include
the first entry but are distinct from those just enumerated is $f(W_{n-2})$.
If $f(W_n) = f(W_{n-1}) + f(W_{n-2})$ then the sequence $W_1, W_2, \ldots$
is at least a $\it{Fibonacci}$ sequence, where each point is the sum
of the two previous points. In fact, this sequence has a rate of growth
greater than the $\it{golden\ ratio}\  {1 + \sqrt{5} \over 2}$, the
(eventual) rate of growth of $\it{Fibonacci}$ sequences.

We examined properties of all permutations up to length 10 and many beyond.  
A pleasantly surprising phenomenon was that ${d \over dn} h(n)$ 
appears to be a monotonically increasing function.  
The permutation
$$5\ 12\ 2\ 7\ 15\ 10\ 4\ 13\ 8\ 1\ 11\ 6\ 14\ 3\ 9$$
has 16874 distinctly ordered subsequences, 
more than $2^{n-1}$ for $n = 15$.
It seemed evident that Wilf's rate of growth could be improved upon.
In this paper we present a class of permutations $\pi_k$ where $f(\pi_k)$
exceeds ${2^{n-2\sqrt{n}} \over \sqrt{n}}$.
 
\section{The Construction}

\begin{definition}
Let $p$ be a permutation.  We call an entry $p_i$ a {\em descent} 
if $p_i > p_{i+1}$.
\end{definition}

\begin{theorem} \label{maint} 
For $n > 3$, there exists an $n$-permutation containing more than 
${2^{n - 2\sqrt{n}} \over \sqrt{n}}$ distinct patterns

\end{theorem}

\begin{proof} 

Let $\pi_k$ denote the permutation
 $$k\ 2k\ \ldots
\ k^2\ \ \ (k - 1)\ (2k - 1)\ \ldots
\ (k^2 - 1)\ \ \  \ldots \ldots
\ \ \ 1\ (k + 1)\ (2k + 1)\ \ldots
\ (k^2 - k + 1)\ \in S_{k^2}$$
For example,

$\pi_3 = 3\ 6\ 9\ \ \ 2\ 5\ 8\ \ \ 1\ 4\ 7$

$\pi_4 = 4\ 8\ 12\ 16\ \ \ 3\ 7\ 11\ 15\ \ \ 2\ 6\ 10\ 14\ \ \ 1\ 5\ 9\ 13$

It should be noted that the only descents in such a permutation are at the last 
entry of each segment, descending to the first entry of the subsequent segment.  As these points play a signifanct role in our proof, we shall denote the first entry of each segment the $\it{perigee}$ of that segment.

Also, each segment is structured so that the $i^{th}$ entry of that segment is less than
the $i^{th}$ entry of each preceding segment.

As counting all patterns of such a permutation leads
to overwhelming complexity, we will restrict our attention to counting only certain patterns.
Let $k > 3$.  For the subsequences under consideration, we require that 
the first $k$ entries of $\pi_k$ be included, i.e., 
every entry in the first segment. 
Also, we include the perigee of each of the other segments.
For example, underlining the entries required in $\pi_4$,
$$\underline{4} \ \underline{8}\ \underline{12}\ \underline{16}\ \ \ 
\underline{3}\ 7\ 11\ 15\ \ \ \underline{2}\ 6\ 10\ 14\ \ \ \underline{1}\ 5\ 9\ 13$$

Requiring these $2k - 1$ entries leaves 
$k^2 - (2k - 1) = (k - 1)^2$ entries to choose from.  
To maximize the number of choices, we will choose half of these 
remaining entries for inclusion.  
Altogether, there will be $\pmatrix{(k - 1)^2 \cr {(k - 1)^2 \over 2}\cr}\ $ 
subsequences of length $(2k - 1) + {(k - 1)^2 \over 2} = {k^2 + 2k - 1 \over 2}$ fulfilling these requirements.
We claim that these subsequences are distinctly ordered, i.e., correspond to distinct patterns.

Suppose $q$ and $r$ are identically ordered subsequences of this type.
Then, the descents in $q$ and $r$ must be at the same positions.
We required the perigee of each segment of $\pi_k$ to be included in $q$ and $r$, 
so any descent will immediately precede a perigee.
Therefore, each perigee occupies the same position in $q$ as in $r$.  This, in turn, implies that the $i^{th}$ entry of $q$ lives in the same segment of $\pi_k$ as does the $i^{th}$ entry of $r$, since they are situated between the same perigees.

Furthermore, included in both $q$ and $r$ are all entries of the first segment of 
$\pi_k$.  
Since $q$ and $r$ are identically ordered, by definition, $q_i < q_j \iff r_i < r_j$.  If there was some entry in the first segment of $\pi_k$ that was less than some $q_i$ and greater than $r_i$, then $q$ and $r$ would not be distinctly ordered as was assumed. 
As noted earlier, the $i^{th}$ entry of each segment of $\pi_k$ is less than the $i^{th}$ entry of $\pi_k$ itself.
So, any two entries $q_i$ and $r_i$ occupy the same position within the same segment and are, in fact, equal. Therefore, $q = r$.  

We've shown that two identically ordered subsequences must actually be the same, and our claim follows that the subsequences are all distinct. 

Stirling's approximation gives
$$\pmatrix{2n \cr n \cr} \rightarrow 
{2^{2n} \over \sqrt{\pi n}}$$

So, the number of distinct subsequences fullfilling the above requirements 
in $\pi_k$ will approach 
$$\sqrt{2 \over \pi} \cdot {2^{(k - 1)^2} \over (k - 1)}$$.

Letting $n = k^2$, the length of the permutation $\pi_k$, 
we have, as $n \rightarrow \infty$,
$$\rm number\ of\ distinct\ patterns\ in \mit\ \pi_k \rightarrow 
\left( {2 \sqrt{2} \over \sqrt{\pi}} \right) \cdot
\left( {2^{n -2\sqrt{n}} \over \sqrt{n} - 1} \right) > 
{2^{n -2\sqrt{n}} \over \sqrt{n}} $$
\end{proof}

\newpage

\section{Long Distance Relationships}
$\indent$
There is still much in this area to be explored.  
While the above class of permutations lends itself to proof,
like Wilf's, it is a tradeoff between manageability and performance.
We have only counted a restricted number of patterns in certain
less than optimal permutations.  The Holy Grail here would be
tighter bounds for $h(n)$.

Second to that a strong result would be an inductive proof on $n$ that
for any permutation $\pi$ of length $n$, one could find a permutation of 
length $n+1$ that contains $\pi$ as well as at least $2 f(\pi)$ 
patterns. This would at least prove that 
$h(n)$ grows at least as fast as $2^n$ for all $n$.

The aim for optimizing the permutation is typically to maximize the
sum of the geographical and numerical distances between any two entries.
That was how the above class of permutations was discovered, although
this property was not explicitly used in the proof.

\end{document}